\documentclass[11pt,reqno]{amsart}

\theoremstyle{plain}
\newtheorem{theorem} {Theorem}[section]
\newtheorem{lemma}[theorem] {Lemma}

\newtheorem{corollary}[theorem] {Corollary}

\theoremstyle{definition}
\newtheorem{definition}[theorem] {Definition}

\theoremstyle{remark}
\newtheorem{remark}[theorem] {Remark}

\numberwithin{equation}{section}

\usepackage{amsfonts}
\usepackage{amssymb}
\usepackage{amscd}

\newcommand{\R}{{\mathbb R}}
\newcommand{\Z}{{\mathbb Z}}
\newcommand{\N}{{\mathbb N}}

\newcommand{\PP}{{\mathcal P}}
\newcommand{\MM}{{\mathcal M}}

\newcommand{\CC}{{\mathbb C}}

\newcommand{\TT}{{\mathcal T}}

\newcommand{\BB}{{\mathfrak B}}
\newcommand{\al}{{\alpha}}
\newcommand{\la}{{\lambda}}
\newcommand{\sa}{{\sigma}}

\newcommand{\iy}{{\infty}}
\newcommand{\vphi}{{\varphi}}
\newcommand{\vep}{{\varepsilon}}
\newcommand{\g}{{\gamma}}
\newcommand{\de}{{\delta}}

\newcommand{\be}{{\beta}}

\newcommand{\bna}{\begin{eqnarray}}
\newcommand{\ena}{\end{eqnarray}}
\newcommand{\ba}{\begin{eqnarray*}}
\newcommand{\ea}{\end{eqnarray*}}
\newcommand{\beq}{\begin{equation}}
\newcommand{\eeq}{\end{equation}}
\DeclareMathOperator*{\esssup}{ess\,sup}

\begin{document}

\title[Constants in Multivariate Inequalities]
{Sharp Constants of Approximation Theory. V. An Asymptotic Equality
Related to Polynomials with Given Newton Polyhedra}
\author{Michael I. Ganzburg}
 \address{Department of Mathematics\\ Hampton University\\ Hampton,
 VA 23668\\USA}
 \email{michael.ganzburg@gmail.com}

 \begin{abstract}
 Let $V\subset\R^m$ be a convex body,
  symmetric about all coordinate hyperplanes, and let $\PP_{aV},\,
  a\ge 0$, be a set of all algebraic polynomials
  whose Newton polyhedra are subsets of $aV$.
 We prove
a limit equality as $a\to \iy$ between the sharp constant in the
multivariate
Markov-Bernstein-Nikolskii
type inequalities
 for polynomials from $\PP_{aV}$ and the corresponding constant for
 entire functions of exponential type with the spectrum in $V$.
 \end{abstract}
 \maketitle

 \section{Introduction}\label{S1}
\setcounter{equation}{0}
\noindent
We continue the study of the sharp constants in multivariate inequalities
of approximation theory
that began in \cite{G2018, G2019, G2019b, G2020}. In this paper we prove
an asymptotic equality
between the sharp constants in the
multivariate
 Markov-Bernstein-Nikolskii
type inequalities for entire functions of exponential type and
algebraic polynomials whose Newton polyhedra are subsets of the given
 convex body.
\vspace{.12in}\\
\textbf{Notation.}
Let $\R^m$ be the Euclidean $m$-dimensional space with elements
$x=(x_1,\ldots,x_m),\, y=(y_1,\ldots,y_m),
\,t=(t_1,\ldots,t_m),\,u=(u_1,\ldots,u_m)$,
the inner product $t\cdot x:=\sum_{j=1}^mt_jx_j$,
and the norm $\vert x\vert:=\sqrt{x\cdot x}$.
Next, $\CC^m:=\R^m+i\R^m$ is the $m$-dimensional complex space with elements
$z=(z_1,\ldots, z_m)=x+iy$
and the norm $\vert z\vert:=\sqrt{\vert x\vert^2+\vert y\vert^2}$;
$\Z^m$ denotes the set of all integral lattice points in $\R^m$;
and $\Z^m_+$ is a subset of $\Z^m$
of all points with nonnegative coordinates.
We also use multi-indices $s=(s_1,\ldots,s_m)\in \Z^m_+,\,
\be=(\be_1,\ldots,\be_m)\in \Z^m_+$,
and $\al=(\al_1,\ldots,\al_m)\in \Z^m_+$
with
 \ba
 \vert s\vert:=\sum_{j=1}^m s_j,\quad
 \vert\be\vert:=\sum_{j=1}^m\be_j,\quad
 \vert\al\vert:=\sum_{j=1}^m\al_j,\quad
 y^\be:=y_1^{\be_1}\cdot\cdot\cdot y_m^{\be_m}, \quad
 D^\al:=\frac{\partial^{\al_1}}{\partial y_1^{\al_1}}\cdot\cdot\cdot
 \frac{\partial^{\al_m}}{\partial y_m^{\al_m}}.
 \ea
 \noindent
Given $\sa\in\R^m,\,\sa_j\ne 0,\,1\le j\le m$, and $M>0$, let
$\Pi^m(\sa):=\{t\in\R^m: \vert t_j\vert\le \vert\sa_j\vert, 1\le j\le m\},\,
Q^m(M):=\{t\in\R^m: \vert t_j\vert\le M, 1\le j\le m\},\,
\BB^m(M):=\{t\in\R^m: \vert t\vert\le M\}$,
and $O^m(M):=\{t\in\R^m: \sum_{j=1}^m\vert t_j\vert\le M\}$
be the $m$-dimensional parallelepiped, cube, ball,
and octahedron, respectively.
In addition, $\vert \Omega\vert_k$ denotes the $k$-dimensional Lebesgue measure
of a  measurable set $\Omega\subseteq\R^m,\,1\le k\le m$.
We also use the floor function
 $\lfloor a \rfloor$.

Let $L_r(\Omega)$ be the space of all measurable complex-valued functions $F$
 on a measurable set $\Omega\subseteq\R^m$  with the finite quasinorm
 \ba
 \|F\|_{L_r(\Omega)}:=\left\{\begin{array}{ll}
 \left(\int_\Omega\vert F(x)\vert^r\, dx\right)^{1/r}, & 0<r<\iy,\\
 \mbox{ess} \sup_{x\in \Omega} \vert F(x)\vert, &r=\iy.
 \end{array}\right.
 \ea
 This quasinorm allows the following "triangle" inequality:
 \beq\label{E1.1}
 \left\|\sum_{j=1}^l F_j\right\|^{\tilde{r}}_{L_r(\Omega)}
 \le \sum_{j=1}^l \left\|F_j\right\|^{\tilde{r}}_{L_r(\Omega)},
 \qquad F_j\in L_r(\Omega),\qquad
 1\le j\le l,
 \eeq
 where $l\in\N:=\{1,\,2,\ldots\}$ and $\tilde{r}:=\min\{1,r\}$ for $r\in(0,\iy]$.

 In this paper we will need certain definitions and properties of
 convex bodies in $\R^m$.
Throughout the paper $V$ is a centrally symmetric (with respect to the origin)
 closed
 convex body in $\R^m$ and
 $V^*:=\{y\in\R^m: \forall\, t\in V, \vert t\cdot y\vert \le 1\}$
 is the \emph{polar} of $V$.
 It is well known that $V^*$ is a centrally symmetric (with respect to the origin)
 closed
 convex body in $\R^m$ and $V^{**} =V$ (see, e.g., \cite[Sect. 14]{R1970}).
 The set $V$ generates the following dual norm
 on $\CC^m$ by
 \ba
  \|z\|_V^*:=\sup_{t\in V}\left\vert\sum_{j=1}^m t_jz_j\right\vert,\quad z\in\CC^m.
 \ea

 Throughout the paper we assume that the body $V\subset \R^m$ satisfies
 the \emph{parallelepiped condition ($\Pi$-condition)}, that is,
 for every vector $t\in V$
 with nonzero coordinates, the parallepiped $\Pi^m(t)$ is a subset of $V$.
 It is easy to verify that $V$ satisfies the $\Pi$-condition if and only if
 $V$ is symmetric about all coordinate hyperplanes, that is,
 for every $t\in V$ the vectors $(\pm\vert t_1\vert,\ldots, \pm\vert t_m\vert)$
 belong to $V$.
   In particular, given $\la\in [1,\iy]$ and
   $\sa\in\R^m,\,\sa_j>0,\,1\le j\le m$,
   the set $V_{\la,\sa}:=\left\{t\in\R^m:
 \left(\sum_{j=1}^m\vert t_j/\sa_j\vert^{\la}\right)^{1/\la}\le 1\right\}$,
 satisfies the $\Pi$-condition.
 Therefore, the sets $\Pi^m(\sa)$ (for $\la=\iy$),
 $Q^m(M)$ (for $\la=\iy$ and $\sa=(M,\ldots,M)$),
 $\BB^m(M)$ (for $\la=2$ and $\sa=(M,\ldots,M)$), and
 $O^m(M)$ (for $\la=1$ and $\sa=(M,\ldots,M)$)
 satisfy the $\Pi$-condition as well.

 Given $a\ge 0$, the set of all trigonometric polynomials
 $T(x)=\sum_{\theta\in aV\cap \Z^m}c_\theta\exp[i(\theta\cdot x)]$
  with complex
 coefficients is denoted by $\TT_{aV}$.

  \begin{definition}\label{D1.1}
 We say that an entire function $f:\CC^m\to \CC^1$ has exponential type $V$
 if for any $\vep>0$ there exists a constant $C_0(\vep,f)>0$ such that
 for all $z\in \CC^m$,
 $\vert f(z)\vert\le C_0(\vep,f)\exp\left((1+\vep)\|z\|_V^*\right)$.
 \end{definition}
  The class of all entire function of exponential type $V$ is denoted
  by $B_V$.
  In the univariate case we use the notation
  $B_\la:=B_{[-\la,\la]},\,\la>0$.
  Throughout the paper, if no confusion may occur, the same notation is applied to
  $f\in B_V$ and its restriction to $\R^m$ (e.g., in the form
  $f\in  B_V\cap L_p(\R^m))$.
  The class $B_V$ was defined by Stein and Weiss
  \cite[Sect. 3.4]{SW1971}. For $V=\Pi^m(\sa),\,V=Q^m(M),$ and
  $V=\BB^m(M)$, similar
  classes were
  defined by Bernstein \cite{B1948} and Nikolskii
   \cite[Sects. 3.1, 3.2.6]{N1969}, see also
  \cite[Definition 5.1]{DP2010}.
  Properties of functions from $B_V$ have been investigated in numerous
  publications (see, e.g., \cite{B1948,  N1969, SW1971, NW1978,
  G1982, G1991, G2001} and
  references therein). Some of these properties are presented in Lemma \ref{L2.1}.

  Given $a\ge 0$, let  $\PP_{aV}$ be a set of all
  polynomials $P(x)=\sum_{\be \in aV\cap \Z_+^m}c_\be x^\be$
   in $m$ variables
   with
  complex coefficients
  whose Newton polyhedra are subsets of $aV$.
  In the univariate case we
  use the notation $\PP_a=\PP_{\lfloor a\rfloor}:=\PP_{a[-1,1]}$.
  In the case of $V=O^m(1),\,\PP_{nV}=\PP_{O^m(n)}$ coincides with the set
   of all polynomials in $m$ variables of total degree at most $n,\,n\in\N$.
  It is easy to verify that if $V_1\subseteq V_2$, then
  $B_{V_1}\subseteq B_{V_2}$ and $\PP_{aV_1}\subseteq\PP_{aV_2}$.

  Throughout the paper $C,\,C_1,\,C_2,\ldots,C_{25}$ denote positive constants independent
of essential parameters.
 Occasionally we indicate dependence on certain parameters.
 The same symbol $C$ does not
 necessarily denote the same constant in different occurrences,
 while $C_k,\,1\le k\le 25$,
 denotes the same constant in different occurrences.
  \vspace{.12in}\\
\textbf{ Markov-Bernstein-Nikolskii Type Inequalities.}
Let
$D_N:=\sum_{\vert\al\vert=N}b_\al D^\al$
be a linear differential
operator with constant coefficients $b_\al\in\CC^1,\,\vert\al\vert=N,\,
N\in \Z^1_+$.
We assume that $D_0$
is the corresponding imbedding or identity operator.

 Next, we define  sharp constants in multivariate Markov-Bernstein-Nikolskii type
inequalities for algebraic and trigonometric polynomials
and entire functions of exponential type. Let
\bna
&&M_{p,D_N,n,m,V}:=n^{-N-m/p}
\sup_{P\in\PP_{O^m(n)}\setminus\{0\}}\frac{\vert D_N(P)(0)\vert}
{\|P\|_{L_p(V^*)}},\label{E1.2a}\\
&&\Tilde{M}_{p,D_N,a,m,V}:=a^{-N-m/p}
\sup_{P\in\PP_{aV}\setminus\{0\}}\frac{\vert D_N(P)(0)\vert}
{\|P\|_{L_p(Q^m(1))}},\label{E1.2}\\
&&P_{p,D_N,a,m,V}:=a^{-N-m/p}
\sup_{T\in\TT_{aV}\setminus\{0\}}\frac{\|D_N(T)\|_{L_\iy(Q^m(\pi))}}
{\|T\|_{L_p(Q^m(\pi))}},\nonumber\\
&& E_{p,D_N,m,V}:=
\sup_{f\in (B_{V}\cap L_p(\R^m))\setminus\{0\}}\frac{\|D_N(f)\|_{L_\iy(\R^m)}}
{\|f\|_{L_p(\R^m)}}.\label{E1.3}
\ena
Here, $a>0,\,N\in\Z^1_+,\,n\in\N,\,V\subset\R^m$, and $p\in(0,\iy]$.
In a sense,
$M_{p,D_N,n,m,V}$ and $\Tilde{M}_{p,D_N,n,m,V},\linebreak
n\in\N$,
 are dual sharp constants
 since the domain of integration $V^*$  in \eqref{E1.2a}
 is the polar of the polynomial "degree" $V$ in \eqref{E1.2},
and the domain of integration  $Q^m(1)=(O^m(1))^*$ in \eqref{E1.2}
 is the polar of the polynomial "degree" $O^m(1)$ in \eqref{E1.2a}.
In particular,
$M_{p,D_N,n,m,O^m(1)}=\Tilde{M}_{p,D_N,n,m,O^m(1)},\,n\in\N$.
We show in this paper that the equality can be asymptotically extended to
any $V$, satisfying the $\Pi$-condition.

Newton polyhedra and polynomial classes $\PP_{aV}$ associated with
Newton polyhedra play an important role in algebra, geometry, and analysis
(see, e.g., a survey \cite[Sect. 3]{AVGK1984}).
However, the only sharp estimate for polynomials from $\PP_{aV}$
we know in multivariate approximation theory is
a sharp V. A. Markov-type inequality
for polynomial coefficients with
$a\in\N$ and $V=\Pi^m(\sa),\,\sa_j\in\N,\,1\le j\le m,$ proved by
Bernstein \cite[Theorem 1]{B1948a} (see \eqref{E1.5c} below).
The purpose of this paper is to prove a limit relation between
 $E_{p,D_N,m,V}$ and
 $\Tilde{M}_{p,D_N,a,m,V}$  as $a\to\iy$
 for $V$, satisfying the $\Pi$-condition.

The following limit relation for multivariate trigonometric polynomials
\beq\label{E1.5}
\lim_{a\to\iy}P_{p,D_N,a,m,V}=E_{p,D_N,m,V},\qquad p\in(0,\iy],
\eeq
was proved by the author \cite[Theorem 1.3]{G2018}.
In the univariate case of $V=[-1,1],\,D_N=d^N/dx^N$, and $a\in\N$,
  \eqref{E1.5} was proved by the author and Tikhonov \cite{GT2017}.
In  earlier publications \cite{LL2015a, LL2015b}, Levin and Lubinsky established
versions of \eqref{E1.5} on the unit circle for $N=0$.
Quantitative estimates of the remainder in asymptotic equalities of the
Levin-Lubinsky type were found by Gorbachev and Martyanov \cite{GM2020}.
Certain extensions of the Levin-Lubinsky's results to the $m$-dimensional
unit sphere in $\R^{m+1}$ were recently proved by Dai, Gorbachev, and Tikhonov
\cite{DGT2018}.

The first sharp constant in the inequality for polynomial coefficients
 was found by V. A. Markov \cite{M1892} (see also
 \cite[Eqs. (5.1.4.1)]{MMR1994})  in the form ($n\in\N$)
 \bna\label{E1.5a}
 &&M_{\iy,d^N/dx^N,n,1,[-1,1]}
 =\Tilde{M}_{\iy,d^N/dx^N,n,1,[-1,1]}\nonumber\\
 &&=\mu^N_n:=n^{-N}
 \left\{\begin{array}{ll}
 \left\vert T_{n-1}^{(N)}(0)\right\vert, &n-N\,\mbox{is odd},\\
 \left\vert T_{n}^{(N)}(0)\right\vert, &n-N\,\mbox{is even}
 \end{array}\right.
\nonumber\\
 &&=1+o(1)=(1+o(1))E_{\iy,d^N/dx^N,1,[-1,1]},
 \ena
 as $n\to\iy$, where $T_n\in\PP_n$ is the Chebyshev polynomial of the first kind.
 For $p=2$ Labelle \cite{L1969} proved the equalities ($n\in\N,\,N\le n$)
 \bna\label{E1.5b}
  &&M_{2,d^N/dx^N,n,1,[-1,1]}
 =\Tilde{M}_{2,d^N/dx^N,n,1,[-1,1]}\nonumber\\
 &&=\frac{(2N)!}{2^N N!}\sqrt{N+1/2}\,\,n^{-(N+1/2)}
  \binom{\lfloor(n-N)/2\rfloor+N+1/2}
  {N+1/2}=\frac{1+o(1)}{\sqrt{\pi(2N+1)}}\nonumber\\
  &&=(1+o(1))E_{2,d^N/dx^N,1,[-1,1]},
 \ena
 as $n\to\iy$.

 The following sharp constant in the
  multivariate inequality for polynomial coefficients
 was found in \cite[Theorem 1]{B1948a}:
 \beq\label{E1.5c}
 \Tilde{M}_{\iy,D^\al,a,m,\Pi^m(\sa)}
 =a^{-\vert\al\vert}
 \prod_{j=1}^m
 \lfloor a\sa_j\rfloor^{\al_j}
 \mu^{\al_j}_{\lfloor a\sa_j\rfloor}
 =(1+o(1))\prod_{j=1}^m\sa_j^{\al_j}
 =(1+o(1))E_{\iy,D^\al,m,\Pi^m(\sa)},
 \eeq
 as $a\to\iy$, where
 $\mu^{\al_j}_{\lfloor a\sa_j\rfloor}$,
  is defined in \eqref{E1.5a} and
  $\sa_j>0,\,1\le j\le m$.
  Note that
 \beq\label{E1.5d}
 \Tilde{M}_{\iy,D^\al,a,m,\Pi^m(\sa)}
 \le \prod_{j=1}^m\sa_j^{\al_j},
 \eeq
 which follows from
 the left  equality in \eqref{E1.5c}
 and the corresponding univariate version of \eqref{E1.5d}
 $\mu^N_{n}\le 1$
 (see \cite[Eq. 2.6(9)]{T1963}
 with its proof in \cite[Lemma 2.5]{G2019b}).

 A crude estimate
 \beq\label{E1.5e}
 \vert c_\be\vert
 \le \left(\prod_{j=1}^m \be_j!\right)^{-1}\,(A(V)a/M)^{\vert\be\vert}
 \|P\|_{L_\iy(Q^m(M))},
 \qquad \be\in aV\cap\Z^m_+,
 \eeq
 for coefficients of a polynomial
 $P(x)=\sum_{\be\in aV\cap\Z^m_+}c_\be x^\be$
 from $\PP_{aV}$
 follows immediately from \eqref{E1.5d}
 if we choose a cube $Q^m(A),\,A=A(V)$, such that
 $V\subseteq Q^m(A)$ and use \eqref{E1.5d} for $\Pi^m(\sa)=Q^m(A)$.

 The author \cite[Theorem 1.2]{G2019b} extended
 \eqref{E1.5a} and \eqref{E1.5b}
  to  a general asymptotic relation for a multivariate
 $L_p$-version of the V. A. Markov constant
 for polynomial coefficients
 in the
 following form ($n\in\N,\,p\in(0,\iy]$):
 \beq\label{E1.6}
 \lim_{n\to\iy}M_{p,D_N,n,m,V}=E_{p,D_N,m,V}.
 \eeq
 For  $m=1,\,D_N=d^N/dx^N$, and $V=[-1,1]$
 this equality was proved
 by the author in \cite[Theorem 1.1]{G2017}.
 A special case of \eqref{E1.6}
 for an even $N\in\Z^1_+,\,p\in[1,\iy]$,
  the unit ball $V=\BB^m(1)$, and the operator
$D_N=\Delta^{N/2}$, where $\Delta$ is the Laplace operator,
was obtained by the author in \cite[Corollary 4.4]{G2019}.

 Note that relations \eqref{E1.5} and \eqref{E1.6} are valid
 for any centrally symmetric $V$ (see \cite{G2018, G2019b}).
 Note also that certain properties of the sharp constants in
univariate weighted spaces are discussed by Arestov and Deikalova
\cite{AD2015}.
In addition, note that the Bernstein-Nikolskii sharp constants
$E_{p,D_N,m,V}$ can be easily found only for $p=2$
(see \cite[Eq. (1.6)]{G2018}).

Despite the fact that the constants
$M_{p,D_N,n,m,V}$ and $\Tilde{M}_{p,D_N,a,m,V}$ for $m>1$
are defined differently by \eqref{E1.2a} and \eqref{E1.2},
it turns out that they are asymptotically equal.
In this paper we extend \eqref{E1.5a}, \eqref{E1.5b},
 and \eqref{E1.5c} to
a general asymptotic relation for $\Tilde{M}_{p,D_N,a,m,V}$,
which is similar to \eqref{E1.6}.
\vspace{.12in}\\
 \textbf{Main Results and Remarks.} Recall that
 $V$ is a
 closed convex body in $\R^m$, satisfying the $\Pi$-condition.
 In particular, $V$ is centrally symmetric (with respect to the origin).

 \begin{theorem} \label{T1.2}
 If $N\in\Z^1_+,\,V\subset\R^m$, and $p\in(0,\iy]$, then
  $ \lim_{a\to\iy}\Tilde{M}_{p,D_N,a,m,V}$ exists and
 \beq \label{E1.7}
  \lim_{a\to\iy}\Tilde{M}_{p,D_N,a,m,V}=E_{p,D_N,m,V}.
 \eeq
 In addition, there exists a nontrivial function $f_0\in  B_V\cap L_p(\R^m)$ such that
 \beq \label{E1.8}
\lim_{a\to\iy}\Tilde{M}_{p,D_N,a,m,V}= \|D_N(f_0)\|_{L_\iy(\R^m)}/\|f_0\|_{L_p(\R^m)}.
 \eeq
 \end{theorem}
 \noindent
 The following corollary is a direct consequence of
 relations \eqref{E1.5}, \eqref{E1.6}, and \eqref{E1.7}.
 \begin{corollary} \label{C1.3}
 If $n\in\N,\,N\in\Z^1_+,\,V\subset\R^m$, and $p\in(0,\iy]$,
 then
 \ba
 \lim_{n\to\iy}M_{p,D_N,n,m,V}=\lim_{a\to\iy}\Tilde{M}_{p,D_N,a,m,V}
 =\lim_{a\to\iy}P_{p,D_N,a,m,V}=E_{p,D_N,m,V}.
 \ea
 \end{corollary}

 \begin{remark}\label{R1.4}
Relations \eqref{E1.7} and \eqref{E1.8} show that the function
$f_0\in B_V\cap L_p(\R^m)$
from Theorem \ref{T1.2} is an extremal function for $E_{p,D_N,m,V}$.
 \end{remark}

\begin{remark}\label{R1.5}
In definitions \eqref{E1.2} and \eqref{E1.3} of the sharp constants we
discuss only complex-valued functions $P$ and $f$. We can define similarly
the "real" sharp constants if the suprema in \eqref{E1.2} and \eqref{E1.3}
 are taken over all real-valued functions
on $\R^m$ from $\PP_{aV}\setminus\{0\}$
and $(B_V\cap L_p(\R^m))\setminus\{0\}$, respectively.
It turns out that the "complex" and "real" sharp constants coincide.
 For $m=1$ this fact was proved in \cite[Sect. 1]{G2017} (cf.
\cite[Theorem 1.1]{GT2017} and \cite[Remark 1.5]{G2019b}),
and the case of $m>1$ can be proved similarly.
\end{remark}

\begin{remark}\label{R1.5a}
Answering a referee's question,
we announced in \cite[Remark 1.6]{G2019b} relation \eqref{E1.7}
for $V=Q^m(M)$ and $a\in\N$ with a typo ($a^{-N-m/p}$ was missing).
\end{remark}

\begin{remark}\label{R1.6}
Now and then we call $M_{p,D_N,n,m,V}$ and $\Tilde{M}_{p,D_N,a,m,V}$
the V. A. Markov constants for polynomial coefficients because of relations
\eqref{E1.5a}. However, there are different constants
$\MM_{p,D_N,n,m,V}$ and $\Tilde{\MM}_{p,D_N,a,m,V}$, defined by
\eqref{E1.2a} and \eqref{E1.2}, respectively,
 with $n^{-N-m/p} \vert D_N(P)(0)\vert$ replaced by
the corresponding $L_\iy$-norm.
 They are associated with the name of
V. A. Markov as well
because he \cite{M1892} found the sharp constant for
$m=1,\,p=\iy,\,D_N=d^N/dx^N$, and $V=[-1,1]$. A brief survey on
$\MM_{p,d^N/dx^N,n,1,[-1,1]}=\Tilde{\MM}_{p,d^N/dx^N,n,1,[-1,1]}$
and its asymptotic behaviour
were presented in \cite{G2017}
(see also \cite[Corollary 4.6]{G2019}).
Certain estimates of $\MM_{p,D_0,n,m,V}$ were surveyed in
\cite[Remark 1.8]{G2019b}.
\end{remark}

The proof of Theorem \ref {T1.2} is presented in
Section \ref{S3}.
It follows general ideas developed in \cite[Corollary 7.1]{G2020}.
Section \ref{S2} contains certain properties of functions
from $B_V$ and $\PP_{aV}$.

\section{Properties of Entire Functions and
  Polynomials}\label{S2}
\setcounter{equation}{0}
\noindent
In this section we discuss certain properties of entire functions
 of exponential type and polynomials
that are needed for the proof of Theorem \ref {T1.2}.
We start with three standard properties of multivariate  entire functions
 of exponential type.
 \begin{lemma}\label{L2.1}
 (a) If $f\in B_V$, then there exists $M=M(V)>0$ such that $f\in B_{Q^m(M)}$.\\
 (b) The following crude Bernstein and Nikolskii type inequalities hold true:
 \bna
 &&\left\|D^\al(f)\right\|_{L_{\iy}(\R^m)}
  \le C
  \left\|f\right\|_{L_{\iy}(\R^m)},\quad f\in B_V\cap L_\iy(\R^m),
  \quad \al\in Z^m_+,\label{E2.1}\\
   &&\left\|f\right\|_{L_{\iy}(\R^m)}
  \le C
  \left\|f\right\|_{L_{p}(\R^m)},\quad f\in B_V\cap L_p(\R^m),\quad
 p\in(0,\iy),\label{E2.2}
  \ena
  where $C$ is independent of $f$.\\
  (c) For any sequence $\{f_n\}_{n=1}^\iy,\,
f_n\in B_V\cap L_\iy(\R^m),\,n\in\N,$
with $\sup_{n\in\N}\| f_n\|_{L_\iy(\R^m)}= C$, there exist a subsequence
$\{f_{n_d}\}_{d=1}^\iy$ and a function $f_0\in B_V\cap L_\iy(\R^m)$
such that for every $\al\in\Z^m_+$,
\beq\label{E2.3}
\lim_{d\to\iy} D^\al f_{n_d}=D^\al f_0
\eeq
uniformly on any compact set in $\CC^m$.
  \end{lemma}
  \proof
   Statement (a) follows from the obvious inclusion $V\subseteq Q^m(M)$
  for a certain $M=M(V)>0$
  (cf. \cite[Lemma 2.1 (a)]{G2019b}).
  Inequality \eqref{E2.1} for $V=Q^m(M),\,M>0$, is well known (see, e.g.,
  \cite[Eq. 3.2.2(8)]{N1969}), while for any $V$, \eqref{E2.1}
  follows from statement (a) (cf. \cite[Lemma 2.1(c)]{G2019b}).
  Inequality \eqref{E2.2} was established in \cite[Theorem 5.7]{NW1978}.
  Statement (c) was proved in \cite[Lemma 2.3]{G2018}.
  \hfill $\Box$

  Given $a\ge 0,\,\g>0$, and a univariate
  continuous function $f\in L_\iy(\R^1)$, let
  \beq\label{E2.3a}
  E(f,\PP_{a},L_\iy([-\g,\g])):=\inf_{R\in\PP_a}\|f-R\|_{L_\iy([-\g,\g])}
  =\|f-R_a\|_{L_\iy([-\g,\g])}
  \eeq
  be the error of best approximation of $f$ by polynomials from $\PP_{a}$
  in the norm of $L_\iy([-\g,\g])$.
  Here,
  $R_a(\cdot)=R_a(f,\g,\cdot)\in\PP_a$ is the polynomial of best uniform
   approximation to $f$.
  Some elementary properties of $R_a$ are discussed in the next lemma.

 \begin{lemma}\label{L2.2}
  (a) The following inequality holds true:
  \beq\label{E2.4}
  \|R_a\|_{L_\iy([-\g,\g])}\le 2\|f\|_{L_\iy(\R^1)}.
  \eeq
  (b) If $f_\mu(v):=f(\mu v),\,\mu\ne 0$, then
  $R_a(f_\mu,\g/\vert\mu\vert,v)
  =R_a(f,\g,\mu v),\,v\in[-\g/\vert\mu\vert,\g/\vert\mu\vert]$.\\
  (c) For $a_j\ge 0,\,\g_j>0$,
  and $t\in\R^d$ with $t_j\ne 0,\,1\le j\le d$,
  the following inequality holds true:
  \bna\label{E2.5}
  &&\max_{\vert x_j\vert\le \g_j/\vert t_j\vert, 1\le j\le d}
  \left\vert\prod_{j=1}^df(t_jx_j)
  -\prod_{j=1}^d R_{a_j}(f,
  \g_j,t_jx_j)\right\vert\nonumber\\
 && \le \|f\|_{L_\iy(\R^1)}^{d-1}
  \sum_{j=1}^d
 2^{j-1} E(f(t_j\cdot),\PP_{a_j},
 L_\iy([-\g_j/\vert t_j\vert,\g_j/\vert t_j\vert])).
  \ena
  \end{lemma}
  \proof
Statement (a) follows from the inequalities
\ba
\|R_a\|_{L_\iy([-\g,\g])}\le \|f\|_{L_\iy([-\g,\g])}
+E(f,\PP_{a},L_\iy([-\g,\g]))\le 2\|f\|_{L_\iy([-\g,\g])},
\ea
while statement (b) is an immediate consequence of the
Kolmogorov characterization of an element of
 best approximation
 to a complex-valued function
\cite{K1948} (see also \cite[Theorem 1.9]{S1974}
 and \cite[Sect. 47]{A1965}).

To prove statement (c), we note that
for $\vert x_j\vert\le \g_j/\vert t_j\vert,\,
1\le j\le d$,
the following relations hold true by \eqref{E2.4}:
\bna
&&\left\vert\prod_{j=1}^df(t_jx_j)
-\prod_{j=1}^d R_{a_j}(f,\g_j,t_jx_j)\right\vert\nonumber\\
&=&\left\vert\sum_{j=1}^d\left[f(t_jx_j)-R_{a_j}
(f,\g_j,t_jx_j)\right]
\prod_{k=j+1}^df(t_kx_k)\prod_{k=1}^{j-1}
R_{a_k}(f,\g_k,t_kx_k)\right\vert
\label{E2.6a}\\
&\le& \|f\|_{L_\iy(\R^1)}^{d-1}
\sum_{j=1}^d
 2^{j-1} \left\|f(t_j\cdot)-R_{a_j}(f,\g_j,t_j\cdot)\right\|_
 {L_\iy([-\g_j/\vert t_j\vert,\g_j/\vert t_j\vert])},\label{E2.6}
\ena
where $\prod_{k=l}^q:=1$ for $q<l$.
Note that the proof of identity \eqref{E2.6a} is simple and
left as an exercise to the reader.
Then \eqref{E2.5} follows from \eqref{E2.6}
since
$R_{a_j}(f,\g_j,t_jx_j)
=R_{a_j}(f_{t_j},\g_j/\vert t_j\vert,x_j)$
by statement (b) for $a=a_j,\,\mu=t_j\ne 0$, and $v=x_j,
\,1\le j\le d$.
\hfill $\Box$

\begin{remark}\label{R2.2a}
Concerning Lemma \ref{L2.2} (b), we note that
for every fixed $v\in\R^1$ the polynomial
$\left\{\begin{array}{ll}
R_a(f,\g,\mu v), &\mu\ne 0,\\
f(0), &\mu=0,
\end{array}\right.$
is obviously a continuous function of
$\mu\in\R^1\setminus \{0\}$,
but it can be discontinuous at $\mu=0$ since
$R_a(f,\g,0)$ is not necessarily equal to $f(0)$.
\end{remark}

  In the next four lemmas we discuss estimates of the error of
  polynomial approximation for functions from $B_V$.

\begin{lemma}\label{L2.3}
Let $g\in B_\la\cap L_\iy(\R^1)$ be a univariate entire function of
exponential type at most $\la>0$.
Given $a\ge 1$ and $\tau\in(0,1)$, the following inequality holds true:
  \beq\label{E2.7}
  E(g,\PP_{a},L_\iy([-a\tau/\la,a\tau/\la]))
  \le C_1(\tau)\exp[-C_2(\tau)\,a]\,\|g\|_{L_\iy(\R^1)},
  \eeq
  where
  \beq\label{E2.8}
  C_1(\tau):=2\left(1+1/\sqrt{1-\tau^2}\right),\qquad
  C_2(\tau):=\log\left(1+\sqrt{1-\tau^2}\right)-\log \tau-\sqrt{1-\tau^2}>0.
  \eeq
  \end{lemma}
  \noindent
  \proof
  It is known (see, e.g., \cite[Sect 5.4.4]{T1963})
  that for any $g\in B_\la\cap L_\iy(\R^1),\,a\ge 1,\,
  \tau\in(0,1)$, and $\de>0$,
  \ba
   E(g,\PP_{a},L_\iy([-a\tau/\la,a\tau/\la]))
   \le \frac{2\exp[a\tau \de]}
   {\de\left(\de+\sqrt{1+\de^2}\right)^{\lfloor a\rfloor}}
   \|g\|_{L_\iy(\R^1)}.
   \ea
   Therefore,
   \bna\label{E2.9}
   &&E(g,\PP_{a},L_\iy([-a\tau/\la,a\tau/\la]))\nonumber\\
   &&\le \frac{2\left(\de+\sqrt{1+\de^2}\right)}{\de}
   \exp\left[\left(\tau\de
   -\log\left(\de+\sqrt{1+\de^2}\right)\right)a\right]
   \|g\|_{L_\iy(\R^1)}.
   \ena
   Setting $\de=\sqrt{1-\tau^2}/\tau$ in \eqref{E2.9},
   we arrive at \eqref{E2.7} and \eqref{E2.8}.
   \hfill $\Box$

   In case of $a\in\N$, versions of  Lemma \ref{L2.3}
  were proved by the author \cite[Lemma 4.1]{G1982}
   and Bernstein \cite[Theorem VI]{B1946}
  (see also \cite[Sect. 5.4.4]{T1963} and
  \cite[Appendix, Sect. 83]{A1965}).  More general and
  more precise inequalities were obtained
  in \cite{G1982} and \cite{G1991}.

  \begin{lemma}\label{L2.4}
  For given $a\ge 1$ and $\tau\in(0,1)$  and for every $t\in V$,
   there exists a polynomial
    $P_t(x)=P_{t,a,V,\tau}(x)
    =\sum_{\be\in aV\cap\Z^m_+}c_{\be}(t)x^\be$
    from $\PP_{aV}$ such that
    $c_{\be}=c_{\be,a,V,\tau}\in L_\iy(V),
    \be\in aV\cap\Z^m_+$, and
     the following inequality holds true:
   \beq\label{E2.10}
   \esssup_{t\in V}\max_{x\in Q^m(a\tau)}
   \vert\exp[i(t\cdot x)]-P_{t}(x)\vert
   \le C_3(\tau,m)\exp[-C_4(\tau,V)\, a],
   \qquad t\in V.
   \eeq
  \end{lemma}
  \proof
  We prove the lemma in three steps.\vspace{.12in}\\
  \textbf{Step 1.}
  We first obtain the univariate inequality ($\la\ne 0$)
  \beq\label{E2.11}
  E\left(\exp[i\la\cdot],\PP_{a},L_\iy
  ([-a\tau/\vert \la\vert,a\tau/\vert\la\vert])\right)
  \le C_1(\tau)\exp[-C_2(\tau)\,a]
  \eeq
  by using Lemma \ref{L2.3} for
  $g(\cdot)=\exp[i\la\cdot]\in B_\la\cap L_\iy(\R^1)$.\vspace{.12in}\\
  \textbf{Step 2.}
Next, we prove \eqref{E2.10} for a parallelepiped
  $V=\Pi^m(u)$, where $u\in\R^m,\,u_j\ne 0,\,1\le j\le m$,
  in the following form:
  \beq\label{E2.12}
  \esssup_{t\in \Pi^m(u)}
  \max_{x\in Q^m(a\tau)}
   \vert\exp[i(t\cdot x)]-P_{t,a,\Pi^m(u),\tau}(x)\vert
   \le C_3(\tau,m)\exp
   \left[-C_4\left(\tau,\Pi^m(u)\right)\, a\right].
   \eeq
   Here,
   \beq\label{E2.13}
   C_3(\tau,m)=m2^{m-1}C_1(\tau),\qquad
   C_4\left(\tau,\Pi^m(u)\right)
   =\min_{1\le j\le m}\vert u_j\vert\,C_2(\tau),
   \eeq
   and the constants $C_1$ and $C_2$ in \eqref{E2.11}
   and \eqref{E2.13}
   are defined by \eqref{E2.8}.
  To prove \eqref{E2.12}, for any $t\in \Pi^m(u)$ we define a
  polynomial
  \beq\label{E2.13a}
  P_t(x)=
  P_{t,a,\Pi^m(u),\tau}(x):=\left\{\begin{array}{ll}
  \prod_{t_j\ne 0,1\le j\le m}
  R_{a\vert u_j\vert}
  \left(\exp[i\cdot],a\tau\vert u_j\vert,t_jx_j\right),
  &t\ne 0,\\
  1,&t=0,
  \end{array}\right.
  \eeq
  from the class $\PP_{a\Pi^m(u)}=\PP_{\Pi^m(au)}$.
  We recall that $R_a=R_a(f,\g,\cdot)$ is defined by \eqref{E2.3a}.
  Since $\vert t_j\vert\le \vert u_j\vert,\,1\le j\le m$,
  we obtain from \eqref{E2.13a},
  \eqref{E2.5}, and \eqref{E2.11}
  \ba
  &&\max_{x\in Q^m(a\tau)}
   \vert\exp[i(t\cdot x)]-P_{t}(x)\vert\\
   &&\le \max_{t_j\ne 0,\,\vert x_j\vert
   \le a\tau\vert u_j\vert/\vert t_j\vert,\,1\le j\le m }
   \vert\exp[i(t\cdot x)]-P_{t}(x)\vert\\
   &&\le 2^{m-1}\sum_{t_j\ne 0,1\le j\le m}
   E(\exp[it_j\cdot],\PP_{a\vert u_j\vert},
   L_\iy\left(\left[-a\tau\vert u_j\vert/
   \vert t_j\vert,a\tau\vert u_j\vert/\vert t_j\vert\right]\right)
   \\
   &&\le m2^{m-1}C_1(\tau)\exp
   \left[-\min_{1\le j\le m}\vert u_j\vert\, C_2(\tau)\, a\right].
   \ea
   This proves \eqref{E2.12} and \eqref{E2.13}.

   Note that by formula \eqref{E2.13a},
   all coefficients $c_\be(t),\,\be\in \Pi^m(au)\cap\Z^m_+$,
   of the polynomial $P_t$  are continuous
   in $t\in \R^m\setminus \bigcup_{j=1}^m H_j$,
   where $H_j$
   is the $j$th $(m-1)$-dimensional coordinate hyperplane
   in $\R^m,\,1\le j\le m$. We also note that the
   coefficients can be discontinuous
   on $H:=\bigcup_{j=1}^m H_j$ (see Remark \ref{R2.2a}).
   However,
   $c_{\be}=c_{\be,a,\Pi^m(au),\tau}\in L_\iy(\R^m),\,
    \be\in \Pi^m(au)\cap\Z^m_+$.
  Indeed,
  using relations \eqref{E2.13a} and \eqref{E2.4},
   we obtain the inequality
   \beq\label{E2.13b}
   \max_{x\in Q^m(a\tau)}\vert P_t(x)\vert
   \le 2^{m}
   \eeq
   for every $t\in \R^m$.
   Therefore, for coefficients of $P_t$ we have the estimate
   $\sup_{t\in \R^m}\vert c_{\be}(t)\vert<\iy,\,
   \be\in \Pi^m(au)\cap\Z^m_+$, by \eqref{E1.5e}
   and \eqref{E2.13b}. Then
   $c_{\be}\in L_\iy(\R^m),\, \be\in \Pi^m(au)\cap\Z^m_+$,
   since $\vert H\vert_m=0$.
    \vspace{.12in}\\
  \textbf{Step 3.}
Finally, let $V$ be a convex body, satisfying the $\Pi$-condition.\\
\textbf{Step 3a).}
   First of all, given $\de\in(1,\iy)$, we construct a finite family of
   parallelepipeds
   \linebreak $\left\{\Pi^m\left(u^{(k)}\right)\right\}_{k=1}^K$
   such that
   \beq\label{E2.14}
   V\subseteq \bigcup_{k=1}^K \Pi^m\left(u^{(k)}\right)\subseteq \de V,\qquad
   \min_{1\le j\le m,1\le k\le K}\left\vert u^{(k)}_j\right\vert
   \ge C_5(\de,V),
   \eeq
   where $K=K(\de,m,V)$.

   To construct the family, we first consider the following
    parallelepipeds
   \ba
   \Pi_l:=\left\{x\in\R^m:\vert x_l\vert\le C_6(\de,V),\,
   \vert x_j\vert\le C_7(\de,V),\,j\ne l\right\},
   \ea
   where $C_6(\de,V):=\min_{1\le l\le m}
   \sqrt{(1+\de)/2}\,\left|OX_l\cap V\right|_1$
   ($OX_l$ is the $l$th coordinate axis, $1\le l\le m$),
    and $C_7(\de,V)$
    is chosen such that
   \beq\label{E2.15}
   \Pi_l\subseteq \sqrt{\de}\, V,\quad 1\le l\le m,
   \qquad
   \inf_{x\in V\setminus\bigcup_{l=1}^m\Pi_l}
   \min_{1\le j\le m}\left\vert x_j\right\vert
   \ge C_7(\de,V).
   \eeq
   Since $\sqrt{(1+\de)/2}<\sqrt{\de}$,
   there is a small enough $C_7(\de,V)<C_6(\de,V)$ such that
   \eqref{E2.15} holds true.

   Next, let $x\in V\setminus\bigcup_{l=1}^m\Pi_l$.
   Then $\vert x_j\vert>0,\,1\le j\le m$, by the second
   relation of \eqref{E2.15}, and
   $x$ is an  interior  point of  $\Pi^m(\de x)$.
   In addition, since $V$ satisfies the $\Pi$-condition,
    we see that $\Pi^m(\de x)\subseteq \de V$.

    Furthermore, setting
   \beq\label{E2.16}
   \Pi_x:=\left\{\begin{array}{ll}
 \sqrt{\de}\,\Pi_l,&x\in \Pi_l\cap V,\,1\le l\le m,\\
 \Pi^m(\de x),
 &x\in V\setminus\bigcup_{l=1}^m\Pi_l,
 \end{array}\right.
   \eeq
   for every $x\in V$,
   we see by the construction of $\Pi_x=\Pi^m(u),\,u=u(x)\in \de V,$
   and by relations \eqref{E2.15} and \eqref{E2.16} that
   \beq\label{E2.16a}
   V\subseteq \bigcup_{x\in V}\Pi_x\subseteq \de V,\qquad
   \min_{x\in V}\min_{1\le j\le m}\left\vert u_j(x)\right\vert
   \ge \sqrt{\de}C_7(\de,V).
   \eeq

   To construct the family
   $\left\{\Pi^m\left(u^{(k)}\right)\right\}_{k=1}^K$
   with $K=K(\de,m,V)$,
   we need the following special case of Morse's theorem
   \cite{M1947} (see also \cite[Remark 1.4]{Gu1975}):

   \begin{lemma}\label{L2.5}
   Let for every $x\in V$ there exist a parallelepiped
   $\tilde{\Pi}_x$ (not necessarily centered at the origin),
    satisfying the condition:
   there exists a fixed constant
   $C\ge 1$  independent of $x$,
   and there exist two balls
   $x+\BB^m(r(x))$ and $x+\BB^m(Cr(x))$
    centered at $x$ of radiuses $r(x)$ and $Cr(x)$,
   respectively, such that
   $x+\BB^m(r(x))\subseteq \tilde{\Pi}_x\subseteq x+\BB^m(Cr(x))$.
   Then a family $\left\{\tilde{\Pi}_x\right\}_{x\in V}$ contains a
   subfamily $\pi:=\left\{\tilde{\Pi}_{x(d)}\right\}_{d=1}^\iy$ with the
   following properties:
   \begin{itemize}
   \item[(a)] $V\subseteq \bigcup_{d=1}^\iy\tilde{\Pi}_{x(d)}$;
   \item[(b)] there exist subfamilies
   $\pi_k,\, 1\le k\le K_1(m,C)$,
   of mutually disjoint parallelepipeds
    such that
   $\pi=\cup_{k=1}^{K_1}\pi_k$.
   \end{itemize}
   \end{lemma}
   \noindent
   Then the family $
   \left\{\tilde{\Pi}_x\right\}_{x\in V}=\left\{\Pi_x\right\}_{x\in V}$ defined by
   \eqref{E2.16} satisfies the condition of Lemma \ref{L2.5}.
   Indeed, by the construction of $\Pi_x$, the condition of
   Lemma \ref{L2.5} is satisfied for
   \ba
   r(x)=C_8(\de,V):=\left(\sqrt{\de}-1\right)C_7,
   \quad x\in V;\qquad
   C=\de D(V)/C_8,
   \ea
   where $D(V)$ is the diameter of $V$.
   In addition, note that any two
   parallelepipeds defined by \eqref{E2.16} have nonempty intersection.
   Hence subfamilies $\pi_k,\, 1\le k\le K_1(m,C)$,
   from property (b) of Lemma \ref{L2.5} contain
   no more than one parallelepiped.
   Therefore, by Lemma \ref{L2.5},
    there exists a finite subfamily of
   parallelepipeds $\left\{\tilde{\Pi}_{x(d)}\right\}_{d=1}^{K_1}
   =\left\{\Pi^m\left(u^{(k)}\right)\right\}_{k=1}^K$
   with $K(\de,m,V):=K_1(m,C)$
   and $u^{(k)}\in \de V,\,1\le k\le K$, such that \eqref{E2.14}
   holds true for $C_5=\sqrt{\de} C_7$
   by \eqref{E2.16a}.\vspace{.12in}\\
  \textbf{Step 3b).}
  Furthermore, given $\tau\in(0,1)$, let us set
   $\de=1/\tau$, and let
  $\left\{\Pi^m\left(u^{(k)}\right)\right\}_{k=1}^K$
  be a finite family of parallelepipeds, where
  $K=K(\de,m,V)$
   and $u^{(k)}\in \de V,\,1\le k\le K$, such that \eqref{E2.14}
   holds true.
   Let us define $P_t(x)=P_{t,a,V,\tau}(x),\,x\in Q^m(a\tau),\,
   t\in V$, by the formula
   \beq\label{E2.17}
   P_{t,a,V,\tau}(x):=P_{t,a/\de,\Pi^m\left(u^{(k)}\right),\tau}(x),
   \qquad  t\in V\cap\left(\Pi^m\left(u^{(k)}\right)
   \setminus \bigcup_{l=1}^{k-1}\Pi^m\left(u^{(l)}\right)\right),
   \quad 1\le k\le K.
   \eeq
   Recall that the polynomial
   $P_{t,a/\de,\Pi^m (u^{(k)}),\tau}(x)$ is defined by \eqref{E2.13a},
   and its coefficients
   belong to
   $L_\iy(\Pi^m (u^{(k)})),\,1\le k\le K$.
   Since $V\subseteq \bigcup_{k=1}^K\Pi^m\left(u^{(k)}\right)$
   by \eqref{E2.14},
   we see from \eqref{E2.17} that the coefficients
   of $\PP_{t,a,V,\tau}$ belong to
   $L_\iy(V)$.
   Next, since
   $\bigcup_{k=1}^K\Pi^m\left(u^{(k)}\right)\subseteq \de V$ by \eqref{E2.14},
   $P_{t,a,V,\tau}\in \PP_{(a/\de)(\de V)}=\PP_{aV}$ for each fixed $t\in V$.

   Furthermore, we obtain from \eqref{E2.17}, \eqref{E2.12},
   \eqref{E2.13}, and \eqref{E2.14}
   \ba
   && \esssup_{t\in V}
   \max_{x\in Q^m(a\tau)}
   \vert\exp[i(t\cdot x)]-P_{t,a,V,\tau}(x)\vert\\
   &&\le  \max_{1\le k\le K}
   \esssup_{t\in \Pi^m\left(u^{(k)}\right)}
   \max_{x\in Q^m(a\tau)}
   \left\vert\exp[i(t\cdot x)]
   -P_{t,a\tau,\Pi^m\left(u^{(k)}\right),\tau}(x)\right\vert\\
   &&\le C_3(\tau,m)
   \exp\left[-\min_{1\le k\le K}\min_{1\le j\le m}
   \vert u_j^{(k)}\vert\,C_2(\tau)\tau\, a\right]\\
   &&\le C_3(\tau,m)\exp
   \left[-C_5(1/\tau,V)C_2(\tau)\tau\, a\right]\\
   &&=C_3(\tau,m)\exp
   \left[-C_4(\tau,V)\, a\right].
   \ea
   This completes the proof of Lemma \ref{L2.4}.
   \hfill $\Box$

   \begin{lemma}\label{L2.6}
  For any $f\in  B_V\cap L_\iy(\R^m),\,\tau\in(0,1)$, and $a\ge 1$,
  there is a polynomial
  $P_a=P_{a,V,\tau,f}\in\PP_{aV}$
   such that for each $\al\in\Z_+^m$ and
    $r\in(0,\iy]$,
   \beq\label{E2.18}
  \lim_{a\to\iy}\left\|D^\al(f)-D^\al(P_a)\right\|_{L_{r}(Q^m(a\tau))}=0.
  \eeq
  \end{lemma}
  \proof
  We prove the lemma in three steps.\vspace{.12in}\\
  \textbf{Step 1.}
  We first assume that $f\in  B_V\cap L_2(\R^m)$.
  By the Paley-Wiener type theorem \cite[Theorem 4.9]{SW1971},
  there exists $\vphi\in L_2(V)$ such that
  $f(x)=(2\pi)^{-m/2}\int_V \vphi(t)
  \exp[i(t\cdot x)]\,dt,\,x\in\R^m$.

 Let $P_t(x)$ be a polynomial from Lemma \ref{L2.4}.
   Then for $a\ge 1$ the integral ($x\in\R^m$)
  \ba
P_a^*(x)=P_a^*(f,V,\tau,x)
:=(2\pi)^{-m/2}\int_V \vphi(t)P_t(x)\,dt
=(2\pi)^{-m/2}\sum_{\be\in aV\cap\Z^m_+}\int_V
\vphi(t)c_\be(t)\,dt\,x^\be
\ea
exists since $c_\be=c_{\be,a,V,\tau}\in L_\iy(V),\,
\be\in aV\cap\Z^m_+$.
Therefore, $P_a^*\in\PP_{aV}$.
Next, it follows from \eqref{E2.10} that given $\tau\in(0,1)$,
\bna\label{E2.19}
\left\|f-P_a^*\right\|_{L_{\iy}(Q^m(a\tau))}
&\le& (2\pi)^{-m/2}\int_V\vert \vphi(t)\vert\, dt\,
\esssup_{t\in V}
\max_{x\in Q^m(a\tau)}
\vert\exp[i(t\cdot x)]-P_{t}(x)\vert\nonumber\\
&\le&  \vert V\vert_m^{1/2}
 C_3(\tau,m)\exp[-C_4(\tau,V)\, a]\,\left\|f\right\|_{L_{2}(\R^m)},
\ena
where $C_3$ and $C_4$ are the constants from
Lemma \ref{L2.4}.\vspace{.12in}\\
 \textbf{Step 2.}
 Next, let  $f\in  B_V\cap L_\iy(\R^m)$.
 Then given $\tau\in (0,1)$ and $\vep\in(0,(1-\tau)/(2\tau C_9)]$,
 where $C_{9}(m,V):=
 (m+1)\sup_{z\in\CC^m}
 \vert z\vert/\|z\|_{V}^*$, the function
 \ba
 f_1(z):=f(z)\left(\frac{\sin\left[\vep
 \left(\sum_{j=1}^mz_j^2\right)^{1/2}\right]}
 {\vep\left(\sum_{j=1}^mz_j^2\right)^{1/2}}\right)^{m+1}
 \ea
 belongs to $B_{(1+\vep C_9)V}\cap L_2(\R^m)$
 and
 $\left\|f_1\right\|_{L_{2}(\R^m)}
 \le C_{10}(m)\,\vep^{-m/2}\left\|f\right\|_{L_{\iy}(\R^m)}$.
 Replacing now $a$ with $a/(1+\vep C_9)$ and
 $\tau$ with $\tau(1+\vep C_9)\le (1+\tau)/2$
 in \eqref{E2.19},
 we see from  \eqref{E2.19} that there exists a polynomial
 $P_a(\cdot)=P_{a,V,\tau,f,\vep}(\cdot)
 :=P^*_{a/(1+\vep C_9)}(f_1,(1+\vep C_9)V,\cdot)\in\PP_{aV}$,
 where $\vep$ will be chosen later,
 such that
 \bna\label{E2.20}
 &&\left\|f_1-P_a\right\|_{L_{\iy}(Q^m(a\tau))}
 \le \left\|f_1-P_a\right\|_{L_{\iy}(Q^m((a/(1+\vep C_9))(1+\tau)/2)}
 \nonumber\\
 &&\le C_{11}(\tau,m,V)\vep^{-m/2}
 \exp[-C_4((1+\tau)/2,V)\, 2a\tau/(1+\tau)]\,\left\|f\right\|_{L_{\iy}(\R^m)}
 \nonumber\\
 &&=C_{11}(\tau,m,V)\vep^{-m/2}
 \exp[-C_{12}(\tau,V)\, a]\,\left\|f\right\|_{L_{\iy}(\R^m)}.
 \ena
 Furthermore, using an elementary inequality
 $v-\sin v\le v^3/6,\,v\ge 0$, we have
 \bna\label{E2.21}
 \left\|f-f_1\right\|_{L_{\iy}(Q^m(a\tau))}
 &\le& (m+1)\max_{x\in Q^m(a\tau)}\left|1-\frac{\sin(\vep\vert x\vert)}
 {\vep\vert x\vert}\right|\left\|f\right\|_{L_{\iy}(\R^m)}\nonumber\\
 &\le& (1/6)(m+1)m\,\vep^2\,a^2\left\|f\right\|_{L_{\iy}(\R^m)}.
 \ena
 Combining \eqref{E2.20} and \eqref{E2.21}, we obtain
 \bna\label{E2.22}
 \left\|f-P_a\right\|_{L_{\iy}(Q^m(a\tau))}
 \le C_{13}(\tau,m,V)\left(\vep^2\,a^2+\vep^{-m/2}
 \exp[-C_{12}(\tau,V)\, a]\right)
 \left\|f\right\|_{L_{\iy}(\R^m)}.
 \ena
 Finally minimizing the right-hand side of \eqref{E2.22} over all
 $\vep\in (0,(1-\tau)/(2\tau C_9)]$,
 we arrive at the following inequality:
\beq\label{E2.23}
 \left\|f-P_a\right\|_{L_{\iy}(Q^m(a\tau))}
 \le C_{14}(\tau,m,V)a^{\frac{2m}{m+4}}
\exp[-C_{15}(\tau,m,V)\, a]
 \left\|f\right\|_{L_{\iy}(\R^m)},
 \eeq
 where $C_{15}=4C_{12}/(m+4)$.
 Note that if the minimum occurs at $\vep=\vep_0$, then
 $P_a=P_{a,V,\tau,f,\vep_0}$ in \eqref{E2.23}.
 \vspace{.12in}\\
 \textbf{Step 3.}
First of all, for $P_b\in\PP_{bV},\,b\ge 1,\,M>0$,
and $\al\in\Z_+^m$,
   we need the following crude Markov-type inequality:
  \beq\label{E2.24}
\left\|D^\al(P_b)\right\|_{L_\iy(Q^m(M))}
  \le C_{16}(m,V,\vert\al\vert)
  (b^2/M)^{\vert\al\vert}\|P_b\|_{L_\iy(Q^m(M))}.
  \eeq
  To prove \eqref{E2.24}, we note that there exists a constant
  $C_{17}(V)$ such that $P_b$ is a polynomial
  of total degree at most
  $n=\lfloor C_{17}b\rfloor\in\N$
  (that is, $P_b\in \PP_{ O^m(n)}$).
  Then inequality \eqref{E2.24} easily follows from a multivariate
  A. A. Markov-type
  inequality proved by Wilhelmsen \cite[Theorem 3.1]{W1974}.

  Next, let $\{P_{a+k}\}_{k=0}^\iy$ be the sequence of
  polynomials, satisfying inequality \eqref{E2.23}
  with $a$ replaced by $a+k,\,k=0,\,1,\,\ldots$.
  Then the series
  \ba
  \sum_{k=0}^\iy \left(P_{a+k+1}-P_{a+k}\right)
  =\lim_{L\to\iy}\left(P_{a+L+1}-P_{a}\right)
  =\lim_{L\to\iy}\left(P_{a+L+1}-f+f-P_{a}\right)
  \ea
  converges to $f-P_a$ in the metric of $L_\iy(Q^m(a\tau))$
  by \eqref{E2.23}.
  In addition, for any $\al\in \Z^m_+$ we obtain by \eqref{E2.24}
  for $M=a\tau$ and by \eqref{E2.23}
  \bna\label{E2.24a}
  &&\sum_{k=0}^\iy \left\|D^\al(P_{a+k+1}-P_{a+k})
  \right\|_{L_\iy(Q^m(a\tau ))}\nonumber\\
  &&\le C_{16}
  (a\tau)^{-\vert\al\vert} \sum_{k=0}^\iy
  (a+k+1)^{2\vert\al\vert}
      \left\|P_{a+k+1}-P_{a+k}\right\|_{L_\iy(Q^m(a\tau ))}
      \nonumber\\
  &&\le C_{16}
  (a\tau)^{-\vert\al\vert}
  \sum_{k=0}^\iy (a+k+1)^{2\vert\al\vert}
      \left(\left\|f-P_{a+k+1}\right\|_{L_\iy(Q^m((a+k+1)\tau ))}
 + \left\|f-P_{a+k}\right\|_{L_\iy(Q^m((a+k)\tau ))}\right)
 \nonumber\\
  &&\le 2C_{14}C_{16}
  (a\tau)^{-\vert\al\vert}\exp[-C_{15}\,a]
  \sum_{k=0}^\iy (a+k+1)^{2\vert\al\vert+2m/(m+4)}
  \exp[-C_{15}\,k]\,\|f\|_{L_\iy(\R^m)}\nonumber\\
  &&\le C_{18}(\tau,m,V,\vert\al\vert,r)\,
  a^{\vert\al\vert+2}\exp[-C_{15}\,a]
  \,\|f\|_{L_\iy(\R^m)}.
  \ena
  Hence the series
  $\sum_{k=0}^\iy D^\al\left(P_{a+k+1}-P_{a+k}\right)$
  is uniformly convergent on $Q^m(a\tau)$ by the Weierstrass
  M-test, and this series converges to
  $D^\al\left(f-P_a\right)$ in the metric of $L_\iy(Q^m(a\tau))$
  by the Differentiation Theorem from multivariate calculus.
  It remains to take account of the following inequalities:
  \bna\label{E2.24b}
 && \left\|D^\al(f)-D^\al(P_a)
  \right\|_{L_r(Q^m(a\tau ))}
  \le (2a\tau)^{m/r}
  \left\|D^\al\left(f-P_a\right)
  \right\|_{L_\iy(Q^m(a\tau ))}\nonumber\\
 && \le (2a\tau)^{m/r}
  \sum_{k=0}^\iy \left\|D^\al(P_{a+k+1}-P_{a+k})
  \right\|_{L_\iy(Q^m(a\tau ))}.
  \ena
  Thus \eqref{E2.18} follows from \eqref{E2.24b} and \eqref{E2.24a},
   and the proof of the lemma is completed.
  \hfill $\Box$

  \begin{remark}\label{R2.8a}
  Note that limit relation \eqref{E2.18} holds true for
   the same polynomial $P_a$ and any $\al\in\Z^m_+$ and
   $r\in(0,\iy]$. The proof of this fact in Lemma \ref{L2.6}
   is based on the exponential approximation rate in \eqref{E2.23}.
   \end{remark}

  A certain polynomial estimate
      is discussed in the following lemma.

    \begin{lemma}\label{L2.7}
     Given  $a\ge 1,\,M>0,\,p\in(0,\iy),\,\tau\in(0,1)$,
     and $P\in\PP_{aV}$, the
    following inequality holds true:
    \beq\label{E2.25}
    \|P\|_{L_\iy(Q^m(\tau M))}
    \le C_{19}(\tau,m,V,p) (a/M)^{m/p}
    \|P\|_{L_p(Q^m(M))}.
    \eeq
    \end{lemma}
    \proof
    Inequality \eqref{E2.25} for $V=Q^m(1)$ and $a\in\N$
    follows from a more general inequality proved in
  \cite[Lemma 2.7 (b)]{G2019b}.
  To prove \eqref{E2.25} for any $V$,
   we note that there exists a constant
  $C_{20}(V)$ such that $P$ is a polynomial of degree at most
  $n=\lfloor C_{20}a\rfloor\in\N$ in each variable
  (that is,
  $P\in
  \PP_{Q^m(n)}$). Then \eqref{E2.25}
  follows from \cite[Lemma 2.7 (b)]{G2019b}.
  \hfill $\Box$

 In the next lemma we discuss special properties of polynomials
 from $\PP_{aV}$.

 \begin{lemma}\label{L2.8}
 Given $a\ge 1,\,b\ge 1$, and
 $P(x)=\sum_{\be \in aV\cap \Z_+^m}c_\be x^\be
 \in\PP_{aV}$, let
 \ba
 R_{a,b}(t):=P(b\sin(t_1/b),\ldots, b\sin(t_1/b)),\qquad t\in\R^m,
 \ea
 be a trigonometric polynomial. Then the following statements
  are valid.\\
 (a) $R_{a,b}\in B_{(a/b)V}$.\\
 (b) For $\al\in\Z^m_+$ the following estimate holds true:
 \beq\label{E2.26}
 \left\vert D^\al(R_{a,b})(0)-D^\al(P)(0)\right\vert
 \le C_{21}(m,\al)\max_{0\le s_j\le \al_j,1\le j\le m,s\ne \al}
 \left\vert D^s(P)(0)\right\vert/b.
 \eeq
 \end{lemma}
 \proof
 (a) We see that
\ba\label{E2.27}
R_{a,b}(bt)
=\sum_{\be\in aV\cap\Z^m_+}
b^{\vert\be\vert}c_\be\prod_{j=1}^m\sin^{\be_j}t_j
=\sum_{\be\in aV\cap\Z^m_+}
b^{\vert\be\vert}c_\be
\sum_{\theta\in\Z^m,0\le\vert \theta_j\vert\le\be_j,1\le j\le m}
d_{\theta,\be}\exp[i(\theta\cdot t)].
\ea
 Then $R_{a,b}(b\cdot)\in\TT_{aV}$,
 since $V$ satisfies the $\Pi$-condition,
and therefore, $R_{a,b}(\cdot)\in B_{(a/b)V}$.\vspace{.12in}\\
(b) To prove this statement, we need the identity
\beq\label{E2.28}
 D^\al(R_{a,b})(0)
 =\sum_{s_1=1}^{\al_1}\ldots \sum_{s_m=1}^{\al_m}
 b^{\vert s\vert-\vert \al\vert}D^s(P)(0)
 \prod_{j=1}^mc(s_j,\al_j),
 \eeq
 where
 \beq\label{E2.29}
 c(l,k):=\sum\frac{k!1^{p_1}(-1)^{p_3}\ldots}
 {p_1!(1!)^{p_1}p_3!(3!)^{p_3}\ldots},
 \eeq
 and the sum in \eqref{E2.29} is taken over all
 nonnegative integers $p_1,\,p_3,\ldots$, such that
 $1p_1+3p_3+\ldots=k$ and $p_1+p_3+\ldots=l
 ,\,0\le l\le k$.

 Identity \eqref{E2.28} for $m=1$ follows from
 Fa\`{a} di Bruno's formula for derivatives
of the composite function
$\psi(b\sin(\cdot/b))$
(see for example \cite{R1980} or \cite{CR1996}).
For $m>1$, \eqref{E2.28}
 can be proved by induction in $m$.

Since $c(k,k)=1,\,k\in\N$, by \eqref{E2.29}, estimate \eqref{E2.26}
follows immediately from \eqref{E2.28}. \hfill $\Box$

\section{Proof of Theorem  \ref{T1.2}}\label{S3}
 \noindent
\setcounter{equation}{0}
Throughout the section we use the notation $\tilde{p}=\min\{1,p\},\,
p\in(0,\iy]$, introduced in Section \ref{S1}.
\vspace{.1in}\\
\emph{Proof of Theorem \ref{T1.2}.}
We first prove the inequality
\beq \label{E3.1}
  E_{p,D_N,m,V}\le\liminf_{a\to\iy}\Tilde{M}_{p,D_N,a,m,V},\qquad p\in(0,\iy].
 \eeq
Let $f$ be any function from $B_V\cap L_p(\R^m),\,p\in(0,\iy]$.
Then $f\in B_{Q^m(M)},\,M=M(V)>0$, by Lemma \ref{L2.1} (a);
 hence
$D_N(f)\in B_{Q^m(M)}$ by \cite[Sect. 3.1]{N1969}
(see also \cite[Lemma 2.1 (d)]{G2018}).
In addition,
$f\in L_\iy(\R^m)$ by
Nikolskii's inequality \eqref{E2.2} and
$D_N(f)\in L_p(\R^m)$ by Bernstein's
and Nikolskii's inequalities \eqref{E2.1} and
\eqref{E2.2} and by the "triangle" inequality
\eqref{E1.1}. Therefore,
\beq\label{E3.2}
\lim_{\vert x\vert\to\iy}D_N(f)(x)=0,\qquad p\in(0,\iy).
\eeq
Indeed, since $D_N(f)\in B_{Q^m(M)}\cap L_p(\R^m)$, \eqref{E3.2} is known for $p\in[1,\iy)$
(see, e.g., \cite[Theorem 3.2.5]{N1969}), and for $p\in(0,1)$ it follows from
\eqref{E2.2}, since if $D_N(f)\in L_p(\R^m),\,p\in(0,1)$,
then $D_N(f)\in L_1(\R^m)$.

Let us first prove \eqref{E3.1} for $p\in(0,\iy)$. Then by \eqref{E3.2}, there exists $x_0\in\R^m$
such that $\|D_N(f)\|_{L_\iy(\R^m)}=\left\vert D_N(f)(x_0)\right\vert$.
Without loss of generality we can assume that $x_0=0$. Let $\tau\in(0,1)$ be a fixed number.
Then using polynomials $P_a\in\PP_{aV},\,a\ge 1$, from Lemma \ref{L2.6}, we obtain
for $r=\iy$ by \eqref{E2.18} and \eqref{E1.2},
\bna\label{E3.3}
 &&\|D_N(f)\|_{L_\iy(\R^m)}=\left\vert D_N(f)(0)\right\vert\nonumber\\
 &&\le  \lim_{a\to\iy}\left\vert D_N(f)(0)-D_N(P_a)(0)\right\vert
 +\liminf_{a\to\iy}\left\vert D_N(P_a)(0)\right\vert\nonumber\\
 &&=\liminf_{a\to\iy}\left\vert D_N(P_a)(0)\right\vert
 \le \tau^{-(N+m/p)}\liminf_{a\to\iy}\left(\Tilde{M}_{p,D_N,a,m,V}
 \left\| P_a\right\|_{L_p(Q^m(a\tau))}\right).
 \ena
 Using again Lemma \ref{L2.6} (for $\al=0$ and $r=p$), we have from \eqref{E1.1}
 \bna\label{E3.4}
 \limsup_{a\to\iy}
 \left\| P_a\right\|_{L_p(Q^m(a\tau))}
 \le \lim_{a\to\iy}\left(\|f-P_a\|_{L_p(Q^m(a\tau))}^{\tilde{p}}
 +\|f\|_{L_p(Q^m(a\tau))}^{\tilde{p}}\right)^{1/\tilde{p}}
 =\|f\|_{L_p(\R^m)}.
 \ena
 Combining \eqref{E3.3} with \eqref{E3.4}, and letting $\tau\to 1-$,
 we arrive at \eqref{E3.1} for $p\in(0,\iy)$.

 In the case $p=\iy$, for any $\vep>0$
 there exists $x_0\in\R^m$
such that $\|D_N(f)\|_{L_\iy(\R^m)}<(1+\vep)\left\vert D_N(f)(x_0)\right\vert$.
Without loss of generality we can assume that $x_0=0$. Then similarly to
\eqref{E3.3} and \eqref{E3.4} we can obtain the inequality
\beq\label{E3.5}
\|D_N(f)\|_{L_\iy(\R^m)}
<(1+\vep)\tau^{-N}\liminf_{a\to\iy}\Tilde{M}_{\iy,D_N,a,m,V}
 \|f\|_{L_\iy(\R^m)}.
 \eeq
 Finally letting $\tau\to 1-$ and $\vep\to 0+$ in \eqref{E3.5},
 we arrive at \eqref{E3.1} for $p=\iy$.
  This completes the proof of \eqref{E3.1}.

  Furthermore, we will prove the
 inequality
 \beq\label{E3.6}
\limsup_{a\to\iy}\Tilde{M}_{p,D_N,a,m,V}\le E_{p,D_N,m,V},\qquad p\in(0,\iy],
\eeq
  by constructing a nontrivial function $f_0\in B_V\cap L_p(\R^m)$,
   such that
 \beq \label{E3.7}
  \limsup_{a\to\iy}\Tilde{M}_{p,D_N,a,m,V}
  \le\|D_N(f_0)\|_{L_\iy(\R^m)}/
\|f_0\|_{L_p(\R^m)}
\le E_{p,D_N,m,V}.
 \eeq
 Then inequalities \eqref{E3.1} and \eqref{E3.6} imply  \eqref{E1.7}.
 In addition, $f_0$ is an extremal function in \eqref{E1.7},
that is, \eqref{E1.8} is valid.

It remains to construct a nontrivial function $f_0$,
satisfying \eqref{E3.7}.
We first note that
\beq \label{E3.8}
\inf_{a\ge 1}\Tilde{M}_{p,D_N,a,m,V}\ge C_{22}(p,N,D_N,m,V).
\eeq
This inequality follows immediately from  \eqref{E3.1}.
Let $U_a\in\PP_{aV}$ be a polynomial, satisfying the equality
\beq \label{E3.9}
\Tilde{M}_{p,D_N,a,m,V}=a^{-N-m/p}\left\vert D_N(U_a)(0)\right\vert
/\|U_a\|_{L_p(Q^m(a))},\qquad a\ge 1.
\eeq
The existence of an extremal polynomial $U_a$ in \eqref{E3.9}
can be proved by the standard compactness argument
(see, e.g., \cite[Proof of Theorem 1.5]{GT2017} and
\cite[Proof of Theorem 1.3]{G2018}).
Next, setting $P_a(x):=U_a(x/a)$,
 we have from \eqref{E3.9} that
\beq \label{E3.10}
\Tilde{M}_{p,D_N,a,m,V}=\left\vert D_N(P_a)(0)\right\vert
/\|P_a\|_{L_p(Q^m(a))}=1/\|P_a\|_{L_p(Q^m(a))},
\eeq
since we can assume that
\beq \label{E3.11}
\left\vert D_N(P_a)(0)\right\vert=1.
\eeq
Then it follows from \eqref{E3.10}, \eqref{E3.11}, and \eqref{E3.8}
that
\ba
\|P_a\|_{L_p(Q^m(a))}
=1/\Tilde{M}_{p,D_N,a,m,V}\le 1/C_{22}(p,N,D_N,m,V).
\ea
Hence using Lemma \ref{L2.7} for $M=a$
 and $\tau\in(0,1)$,
we obtain  the estimate
\beq \label{E3.12}
\sup_{a\ge 1}\|P_a\|_{L_\iy(Q^m(a\tau))}
\le C_{19}/C_{22}
= C_{23}(\tau,p,N,D_N,m,V).
\eeq
In addition,  combining estimates \eqref{E1.5e} for $M=a\tau$ and
\eqref{E3.12}, we have for any $s\in\Z^m_+$,
\beq \label{E3.13}
\left\vert D^s(P_a)(0)\right\vert
\le (A(V)/\tau)^{\vert s\vert} C_{23}
= C_{24}(\tau,p,N,D_N,m,V,s).
\eeq
Furthermore, we define a trigonometric polynomial
\ba
R_{a,a\tau}(t)
:=P_a(a\tau\sin(t_1/(a\tau)),\ldots, a\tau\sin(t_m/(a\tau))),\qquad
t\in \R^m.
\ea
Then $R_{a,a\tau}$ satisfies the following properties:
\begin{itemize}
\item[(P1)] $R_{a,a\tau}\in B_{(1/\tau)V}$.
\item[(P2)] The following relations hold true:
\beq \label{E3.14}
\sup_{a\ge 1}\|R_{a,a\tau}\|_{L_\iy(Q^m(a\tau\pi/2)}
=\sup_{a\ge 1}\|R_{a,a\tau}\|_{L_\iy(\R^m)}\le C_{23}.
\eeq
\item[(P3)] For $\al\in\Z^m_+$ and $a\tau\ge 1$,
\bna \label{E3.15}
 \left\vert D^\al(R_{a,a\tau})(0)-D^\al(P_{a})(0)\right\vert
 &\le& C_{21}\max_{0\le s_j\le \al_j,1\le j\le m,s\ne \al}
 C_{24}(\tau,p,N,D_N,m,V,s)/(a\tau)\nonumber\\
 &=&C_{25}(\tau,p,N,D_N,m,V,\al)/a.
\ena
\item[(P4)] For $a\ge 1,\,p\in(0,\iy]$, and
$M\in(0,a\tau/\sqrt{m}]$,
\beq \label{E3.16}
\|P_a\|_{L_p(Q^m(a)}
\ge \left(1-mM^2(a\tau)^{-2}\right)^{1/p}
\|R_{a,a\tau}\|_{L_p(Q^m(M))}.
\eeq
\end{itemize}
Indeed, property (P1) follows from Lemma \ref{L2.8} (a),
while (P2) is an immediate consequence of \eqref{E3.12}.
Next, property (P3) follows from Lemma \ref{L2.8} (b) and
relations \eqref{E3.13}. To prove (P4), we note first that
for $p=\iy$ inequality \eqref{E3.16} is trivial. Next,
setting $f(\cdot)=\cos(\cdot)$ and replacing
$R_{a_j},\,1\le j\le m$, with $1$ in
identity \eqref{E2.6a},
we obtain
\bna\label{E3.17}
1-\prod_{j=1}^m \cos(t_jx_j)
= \left\vert\sum_{j=1}^m\left(1-\cos(t_jx_j)\right)
\prod_{k=j+1}^m\cos(t_kx_k)\right\vert
\le (1/2)\sum_{j=1}^m(t_jx_j)^2.
\ena
Furthermore, for $a\ge 1,\,p\in(0,\iy)$, and
$M\in(0,a\tau/\sqrt{m}]$,
\bna\label{E3.18}
&&\|P_a\|^p_{L_p(Q^m(a))}
\ge \|P_a\|^p_{L_p(Q^m(a\tau))}
=\int_{Q^m(a\tau\pi/2)}\left\vert
R_{a,a\tau}(t)\right\vert^p \prod_{j=1}^m
\cos(t_j/(a\tau))\,dt\nonumber\\
&&\ge \|R_{a,a\tau}\|^p_{L_p(Q^m(M))}
-\int_{Q^m(M)}\left\vert
R_{a,a\tau}(t)\right\vert^p \left(1-\prod_{j=1}^m
\cos(t_j/(a\tau))\right)\,dt.
\ena
Finally using estimate \eqref{E3.17}
for $x_j=1/(a\tau),\,1\le j\le m$, we obtain
\eqref{E3.16} from \eqref{E3.18}.

Let $\{a_n\}_{n=1}^\iy$ be an increasing sequence
of numbers such that
 $\inf_{n\in\N}a_n\ge 1,\,\lim_{n\to\iy}a_n=\iy$,
and
\beq \label{E3.19}
\limsup_{a\to\iy}\Tilde{M}_{p,D_N,a,m,V}
=\lim_{n\to\iy}\Tilde{M}_{p,D_N,a_n,m,V}.
\eeq

Property (P1) and relation \eqref{E3.14}
of property (P2)
 show that the  sequence of trigonometric polynomials
$\{R_{a_n,a_n\tau}\}_{n=1}^\iy
=\{f_n\}_{n=1}^\iy$
 satisfies the conditions
of Lemma \ref{L2.1} (c)
with $B_V$ replaced by $B_{(1/\tau)V}$.
Therefore, there exist
a subsequence
 $\{R_{a_{n_d},a_{n_d}\tau}\}_{d=1}^\iy$ and
 a function $f_{0,\tau}\in B_{(1/\tau)V}$ such that
 \beq \label{E3.20}
  \lim_{d\to\iy}R_{a_{n_d},a_{n_d}\tau}=f_{0,\tau},\qquad
 \lim_{d\to\iy}D_N\left(R_{a_{n_d},a_{n_d}\tau}\right)
 =D_N(f_{0,\tau}),
 \eeq
 uniformly on any cube $Q^m(M),\,M>0$.

 Moreover, by
  \eqref{E3.11}, \eqref{E3.15}, and \eqref{E3.20},
  \beq \label{E3.21}
  \left\vert D_N(f_{0,\tau})(0)\right\vert
  =\lim_{d\to\iy}\left\vert D_N
  \left(R_{a_{n_d},a_{n_d}\tau}\right)(0)
  \right\vert
  =\lim_{d\to\iy}\left\vert D_N
  \left(P_{a_{n_d}}\right)(0)
  \right\vert
  =1.
\eeq
In addition,
using \eqref{E1.1}, \eqref{E3.20}, \eqref{E3.16},
\eqref{E3.10}, and \eqref{E3.19},
we obtain for any cube $Q^m(M),\,M>0$,
\bna \label{E3.22}
&&\|f_{0,\tau}\|_{L_p(Q^m(M))}
\le \lim_{d\to\iy}
\left(\left\|f_{0,\tau}-R_{a_{n_d},a_{n_d}\tau}
\right\|_{L_p(Q^m(M))}^{\tilde{p}}
+\left\|R_{a_{n_d},a_{n_d}\tau}\right\|_{L_p(Q^m(M))}^{\tilde{p}}
\right)^{1/\tilde{p}}\nonumber\\
&&=\lim_{d\to\iy}\left\|R_{a_{n_d},a_{n_d}\tau}\right\|_{L_p(Q^m(M))}
\le \lim_{d\to\iy}\left\|P_{a_{n_d}}\right\|_
{L_p\left(Q^m\left(a_{n_d}\right)\right)}
=1/
\lim_{d\to\iy}\Tilde{M}_{p,D_N,a_{n_d},m,V}.
\ena
Next using \eqref{E3.22} and \eqref{E3.8}, we see that
\beq \label{E3.23}
\|f_{0,\tau}\|_{L_p(\R^m)}\le 1/C_{22}(p,N,D_N,m,V).
\eeq
Therefore, $f_{0,\tau}$ is a nontrivial function from
$B_{(1/\tau)V}\cap L_p(\R^m)$,
by \eqref{E3.23} and \eqref{E3.21}.
Thus for any cube $Q^m(M),\,M>0$,
we obtain from \eqref{E3.19}, \eqref{E3.10},
 \eqref{E3.16}, \eqref{E3.20}, and \eqref{E3.21}
\bna \label{E3.24}
\limsup_{a\to\iy}\Tilde{M}_{p,D_N,a,m,V}
&=&\lim_{d\to\iy}\left(\left\|P_{a_{n_d}}\right\|_
{L_p\left(Q^m\left(a_{n_d}\right)\right)}\right)^{-1}\nonumber\\
&\le& \lim_{d\to\iy}\left(\left\| R_{a_{n_d},a_{n_d}\tau}\|\right\|
_{L_p(Q^m(M))}\right)^{-1}\nonumber\\
&=&\left\vert D_N(f_{0,\tau})(0)\right\vert/
\|f_{0,\tau}\|_{L_p(Q^m(M))}.
\ena
It follows from \eqref{E3.24} that
\beq \label{E3.25}
\limsup_{a\to\iy}\Tilde{M}_{p,D_N,a,m,V}
\le E_{p,D_N,m,(1/\tau)V}=\tau^{-N-m/p}E_{p,D_N,m,V}.
\eeq
Then letting $\tau\to 1-$ in
\eqref{E3.25}, we arrive at \eqref{E3.6}.
 However, we need to prove  stronger relations
\eqref{E3.7}.

To construct $f_0$, note first that
$f_{0,\tau}(\tau\cdot)\in B_V$
and by \eqref{E3.23}
and \eqref{E2.2},
\ba
\sup_{\tau\in (1/2,1)}\|f_{0,\tau}\|_{L_\iy(\R^m)}
=\sup_{\tau\in (1/2,1)}\|f_{0,\tau}(\tau\cdot)\|_{L_\iy(\R^m)}
\le C \sup_{\tau\in (1/2,1)}
\tau^{-m/p}\|f_{0,\tau}\|_{L_p(\R^m)}
< \iy.
\ea
Therefore, by Lemma \ref{L2.1} (c) applied to a sequence
$\left\{f_{0,\tau_n}\right\}_{n=1}^\iy$, where
$\tau_n\in(1/2,1),\,n\in\N$, and $\lim_{n\to\iy}\tau_n=1$,
there exist a subsequence
$\{f_{0,\tau_{n_d}}\}_{d=1}^\iy$ and a function
$f_0\in B_V\cap L_\iy(\R^m)
=\bigcap_{d=1}^\iy \left(B_{(1/\tau_{n_d})V}
\cap L_\iy(\R^m)\right)$
such that for every $\al\in\Z^m_+$,
$\lim_{d\to\iy} D^\al \left(f_{0,\tau_{n_d}}\right)
=D^\al \left(f_0\right)$
uniformly on any compact set in $\CC^m$.

Note that by \eqref{E3.21} and \eqref{E3.23},
$f_0$ is a nontrivial function from
$B_{V}\cap L_p(\R^m)$. Then using \eqref{E3.24},
we obtain
\ba
\limsup_{a\to\iy}\Tilde{M}_{p,D_N,a,m,V}
&\le& \lim_{M\to\iy}\lim_{n\to\iy}
\left\vert D_N(f_{0,\tau_n})(0)\right\vert/
\left\|f_{0,\tau_n}\right\|_{L_p(Q^m(M))}\\
&=&\left\vert D_N(f_{0})(0)\right\vert/
\left\|f_{0}\right\|_{L_p(\R^m)}\\
&\le& \left\| D_N(f_{0})\right\|_{L_\iy(\R^m)}/
\|f_{0}\|_{L_p(\R^m)}\\
&\le& E_{p,D_N,m,V}.
\ea
Thus \eqref{E3.7} holds true, and this completes
the proof of the theorem.
\hfill$\Box$

\begin{remark}\label{R3.1}
The proof of \eqref{E3.1} is all but identical to the proof of the
inequality
$ E_{p,D_N,m,V}
\le\liminf_{n\to\iy}{M}_{p,D_N,n,m,V}$
from \cite{G2019b} though these proofs are based on different lemmas.
However, the proof of \eqref{E3.6} is different compared with
the proof of the
inequality
$ \limsup_{n\to\iy}{M}_{p,D_N,n,m,V}
\le E_{p,D_N,m,V}$ from \cite{G2019b}.
 The latter proof is based on V. A. Markov-type
inequalities for polynomials from $\PP_{O^m(n)}$.
We do not know if there are analogues of these
inequalities for polynomials from  $\PP_{aV}$
but the case of $V=\Pi^m(\sa)$, see \eqref{E1.5c}.
That is why, \eqref{E3.6} is reduced to certain
relations for trigonometric polynomials
(cf. \cite{G2018}).
\end{remark}
\noindent
\textbf{Acknowledgements.} We are grateful to both anonymous referees
 for valuable suggestions.


\begin{thebibliography}{99}

\bibitem{A1965} N. I. Akhiezer, Lectures on the Theory of
Approximation, second ed., Nauka, Moscow, 1965 (in Russian).

 \bibitem{AD2015} V. Arestov, M. Deikalova, Nikol'skii inequality between
 the uniform norm and $L_q$-norm with ultraspherical weight of
 algebraic polynomials on an interval, Comput. Methods Funct. Theory
  \textbf{15}  (2015), 689--708.

  \bibitem{AVGK1984} V. I. Arnold, A. N. Varchenko, A. B. Givental,
  A. G. Khovanskii,  Singularities of functions, wave fronts,
  caustics and multidimensional integrals. Translated from the Russian.
  Soviet Sci. Rev. Sect. C Math. Phys. Rev., 4, Mathematical physics reviews,
  Vol. 4, 1--92, Harwood Academic Publ., Chur, 1984.

\bibitem{B1946} S. N. Bernstein, On the best approximation of continuous
 functions on the whole real axis by entire functions of given degree, V,
 Dokl. Akad. Nauk SSSR \textbf{54} (1946), 479--482 (in Russian).


  \bibitem{B1948a} S. N. Bernstein, On certain elementary extremal
  properties of polynomials in several variables,
 Dokl. Akad. Nauk SSSR \textbf{59} (1948), 833--836 (in Russian).

  \bibitem{B1948} S. N. Bernstein, On entire functions of finite degree of
 several variables,
 Dokl. Akad. Nauk SSSR \textbf{60} (1948), 949--952 (in Russian).

  \bibitem{CR1996} G. M. Constantine, T. H. Savits,
A multivariate Fa\`{a} di Bruno formula with applications, Trans. Amer. Math. Soc.
\textbf{348} (1996), no. 2, 503--520.

 \bibitem{DGT2018} F. Dai, D. Gorbachev, S. Tikhonov, Nikolskii constants for
 polynomials on the unit sphere, J. d'Analyse Math.
 \textbf{140} (2020), 161--185.

 \bibitem{DP2010} Z. Ditzian, A. Prymak, On Nikolskii
 inequalities for Lorentz spaces,
 Rocky Mountain Journal of Mathematics \textbf{40} (2010), no. 1,
209--223.

 \bibitem{G1982} M. I. Ganzburg, Multidimensional limit theorems of the theory
of best polynomial approximations, Sibirsk. Mat. Zh. \textbf{23} (1982),
no. 3, 30--47 (in Russian); English transl. in
 Siberian Math. J. \textbf{23} (1983), no. 3, 316--331.

 \bibitem{G1991} M. I. Ganzburg,	Limit theorems for the best
 polynomial approximation in the $L_\iy$-metric, Ukrain. Mat. Zh.
 \textbf{23} (1991), no. 3,
  336--342 (in Russian); English transl. in Ukrainian Math. J. \textbf{23}
  (1991), 299--305.

  \bibitem{G2001} M. I. Ganzburg, Polynomial inequalities on measurable sets and
their applications, Constr. Approx. \textbf{17} (2001), 275--306.

 \bibitem{G2017} M. I. Ganzburg, Sharp constants in V. A. Markov-Bernstein
  type inequalities of different metrics, J. Approx. Theory
   \textbf{215} (2017), 92--105.

  \bibitem{G2018} M. I. Ganzburg, Sharp constants of approximation theory.
  I. Multivariate Bernstein-Nikolskii type inequalities,
  J. Fourier Anal. Appl. \textbf{26}:11 (2020),
  doi:10.1007/s00041-019-09720-x .

\bibitem{G2019} M. I. Ganzburg, Sharp constants of approximation theory.
  II. Invariance theorems and  certain
  multivariate inequalities of different metrics,
   Constr. Approx. \textbf{50} (2019), 543--577.

 \bibitem{G2019b} M. I. Ganzburg, Sharp constants of approximation theory.
  III. Certain polynomial  inequalities of different metrics on convex sets,
   J. Approx. Theory \textbf{252} (2020), doi:10.1016/j.jat.2019.105351.

    \bibitem{G2020} M. I. Ganzburg, Sharp constants of approximation theory.
  IV. Asymptotic Relations in General Settings,
  submitted; arXiv:2002.10512.

\bibitem{GT2017} M. I. Ganzburg, S. Yu. Tikhonov, On sharp constants in
Bernstein-Nikolskii inequalities, Constr. Approx. \textbf{45} (2017), 449--466.

\bibitem{GM2020} D. V. Gorbachev, I. A. Martyanov, Bounds of the Nikolskii
polynomial constants in $L^p$ with a Gegenbauer weight,
Trudy Inst. Mat. Mekh. UrO RAN, 26 (2020), no. 4,
126--137 (in Russian).

\bibitem{Gu1975} M. de Guzm\'{a}n, Differentiation of Integrals in $R^n$,
Lecture Notes of Mathematics, \textbf{481}, Springer-Verlag, New-York, 1975.

\bibitem{K1948} A. N. Kolmogorov, A remark on the polynomials of P. L. Chebyshev
deviating the least from a given function, Uspehi Matem. Nauk,
 \textbf{3}  (1948), no. 1, 216--221 (in Russian).

\bibitem{L1969} G. Labelle, Concerning polynomials on the unit interval,
  Proc. Amer. Math. Soc. \textbf{20} (1969),  321--326.

 \bibitem{LL2015a} E. Levin, D. Lubinsky, $L_p$ Christoffel functions, $L_p$
 universality, and Paley-Wiener spaces, J. D'Analyse Math. \textbf{125} (2015), 243--283.

 \bibitem{LL2015b} E. Levin, D. Lubinsky, Asymptotic behavior of Nikolskii constants
 for polynomials on the unit circle, Comput. Methods Funct. Theory \textbf{15} (2015),
 459--468.

 \bibitem{M1892} V. A. Markov, On Functions Deviating Least from Zero
  in a Given Interval, Izdat. Imp. Akad. Nauk, St. Petersburg, 1892
  (in Russian); German transl. in Math. Ann. \textbf{77} (1916), 213--258.

 \bibitem{MMR1994} G. V. Milovanovi\'{c}, D. S. Mitrinovi\'{c}, Th. M. Rassias,
 Topics in Polynomials: Extremal Problems, Inequalities, Zeros,
 World Scientific, Singapore, 1994.

 \bibitem{M1947} A. P. Morse, Perfect blankets, Trans. Amer. Math. Soc.,
 \textbf{61} (1947), 418--442.

 \bibitem{NW1978} R. J. Nessel, G. Wilmes, Nikolskii-type inequalities
  for trigonometric polynomials and entire functions of exponential type,
  J. Austral. Math. Soc., Ser. A \textbf{25} (1978), 7--18.

 \bibitem{N1969} S. M. Nikolskii, Approximation of Functions of Several Variables
 and Imbedding Theorems, Nauka, Moscow, 1969 (in Russian); English edition:
 Die Grundlehren der Mathematischen Wissenschaften, Band 205, Springer-Verlag,
 New York-Heidelberg, 1975.

\bibitem{R1970} R. T. Rockafellar, Convex Analysis, Princeton University Press,
Princeton, 1970.

 \bibitem{R1980} S. Roman,  The formula of Faa di Bruno,
 Amer. Math. Monthly \textbf{87} (1980), 805--809.

 \bibitem{S1974} I. Singer, The Theory of Best Approximation and Functional Analysis,
Conference Board of the Mathematical Sciences Regional Conference Series in Applied Mathematics, No. 13,
Society for Industrial and Applied Mathematics, Philadelphia, Pa., 1974.

  \bibitem{SW1971} E. Stein, G. Weiss, Introduction to Fourier Analysis on
  Euclidean Spaces, Princeton Univ. Press,
   Princeton, 1971.

 \bibitem{T1963} A. F. Timan, Theory of Approximation of Functions
  of a Real Variable, Pergamon Press, New York, 1963.

\bibitem{W1974} D. R. Wilhelmsen, A Markov inequality in several dimensions,
J. Approx. Theory \textbf{11} (1974), 216--220.



\end{thebibliography}
\end{document}